\documentclass[12pt,oneside]{amsart}

\usepackage{amssymb}
\usepackage[french]{babel}
\usepackage[latin1]{inputenc}
\usepackage{graphicx}

\def\N{\mathbb{N}}
\def\Z{\mathbb{Z}}
\def\R{\mathbb{R}}

\def\a{\alpha}

\def\g{\gamma}
\def\d{\delta}
\def\e{\varepsilon}

\def\l{\lambda}

\def\w{\omega}
\def\G{\Gamma}

\def\cal{\mathcal}

\def\ov{\overline}
\def\ol{\overline}

\def\lb{\lbrack}
\def\rb{\rbrack}

\def\lgr{\mathrm{lgr}}

\def\el{\'el\'ement}

\newtheorem{thm}{Théorème}[section]
\newtheorem{cor}{Corollaire}[section]
\newtheorem{lem}{Lemme}[section]
\newtheorem{prop}{Proposition}[section]

\title{Solutions à divers problèmes de décision dans un groupe hyperbolique}

\begin{document}
\maketitle 
\begin{center}
{\sc Jean-Philippe PR\' EAUX}\footnote[1]{Centre de recherche de
l'Armée de l'Air, Ecole de l'air, F-13661 Salon de Provence air}\
\footnote[2]{Centre de Math\'ematiques et d'informatique,
Universit\'e de Provence, 39 rue F.Joliot-Curie, F-13453 marseille
cedex 13\\
\indent {\it E-mail :} \ preaux@cmi.univ-mrs.fr}
\end{center}
\begin{abstract}
Nous établissons des solutions algorithmiques à divers problèmes de décision dans un groupe hyperbolique au sens de
Gromov.
\end{abstract}

\section*{Introduction}

Nous établissons divers algorithmes dans un groupe hyperbolique $G$ : déterminer le centre de $G$, ses éléments de
torsion, les racines d'un \el\ d'ordre infini, son centralisateur, le centralisateur d'un sous-groupe de type fini,
maximalité et malnormalité de sous-groupes cycliques infinis, solution au problème du mot généralisé d'un sous-groupe
virtuellement $\Z$, {\it etc...}

 Notre étude est entièrement basée sur une poignée de propriétés
élémentaires de l'hyperbolicité : d'une part deux propriétés concernant les quasigéodésiques dans un espace
hyperbolique (stabilité des quasigéodésiques de longueur infinie et géodésicité locale implique quasigéodésicité
globale) et d'autre part une solution au problème du mot et de la conjugaison dans un groupe hyperbolique.

\section{Rappels sur les groupes hyperboliques}

Il existe de multiples fa\c cons  de d\'efinir la notion de groupe
et d'espace hyperbolique au sens de Gromov (cf. \cite{gro},
\cite{cdp}, \cite{gdlh}, \cite{sho}). Nous utilisons l'approche
caract\'erisant l'hyperbolicit\'e par la propriété dite  des
''triangles fins''.
Nous nous restreindrons dans la suite \`a des groupes de type
fini.
Soit $G$ un groupe de type fini, et $S$ une famille g\'en\'eratrice finie pour $G$. Nous notons $S^*$ le monoïde libre
sur $S$. Un \el\ $\w$ de $S^*$ est appelé un $\w$ un mot sur $S$, \emph{i.e.} $\w\equiv s_1\cdots s_i \cdots s_n$, avec
$n\geq 0$ et pour tout $i=1\cdots n$, $s_i\in S\cup S^{-1}$ ; il représente un \el\ de $G$. Nous notons $\lgr_S (\w)$,
la longueur du mot $\w$, c'est \`a dire l'entier $n$. Soit $g$ un \el\ de $G$. On d\'efinit l'entier $|g|$ par,
$$|g|=inf\{\lgr_S(\w)| \w\;\mathrm{est\; un\; mot\; repr\acute{e}sentant}\, g\}$$
Clairement cette borne est atteinte, c'est \`a dire qu'il existe
un mot $\w$ de longueur $|g|$, repr\'esentant $g$ (en g\'en\'eral
non unique). Nous d\'efinissons alors une distance sur $G$,
appelée {\it métrique du mot}, en posant pour tout $g_1,g_2\in G$,
$$d_S(g_1,g_2)=|g_1^{-1}g_2|$$
Cette m\'etrique d\'epend de la famille g\'en\'eratrice $S$
considérée.
Cependant lorsqu'une famille g\'en\'eratrice $S$ sera suppos\'e
fix\'ee, nous aurons coutume de noter $\lgr$ et $d$ au lieu de
$\lgr_S$ et $d_S$ et commettrons l'abus de langage consistant \`a
noter $|\w|$ pour un mot $\w$. Plus g\'en\'eralement, donn\'e un
mot, nous le confondrons souvent avec l'\el\ qu'il repr\'esente
dans le groupe.\\

 Donné un groupe $G$ un groupe muni d'une famille g\'en\'eratrice finie $S$, on
 d\'efinit un graphe connexe, orient\'e et localement fini,
  appel\'e le \textit{graphe de Cayley}, not\'e $\G (G,S)$, de
  la fa\c con suivante :\\
-- Les sommets de $\G (G,S)$ sont en bijection avec les \el s de
  $G$. Si $g\in G$, on note $\ol{g}$ le sommet de $\G (G,S)$ lui
  correspondant.\\
-- Il existe une ar\^ete ayant pour origine $\ol{g_1}$ et pour extr\'emit\'e
  $\ol{g_2}$, lorsque il existe $s\in S\cup S^{-1}$ tel que $g_2=g_1.s$
 dans
  $G$. On munit alors cette ar\^ete du label $s$.

  Un chemin a naturellement pour label un mot $\w$ sur $S$ obtenu par concat\'enation des labels des
ar\^etes successives le composant. Il existe dans $\G (G,S)$ un
chemin de label $\w$ du sommet $g_1$ au sommet $g_2$ si et
seulement si $g_2=g_1w$ dans $G$.

On munit $\G (G,S)$ de la m\'etrique simpliciale en assignant à
chaque ar\^ete la longueur 1. Cela fait de $\G (G,S)$ un espace
m\'etrique g\'eod\'esique propre. L'ensemble de ses sommets muni
de la m\'etrique induite est isom\'etrique \`a $(G,d_S)$. Si un
chemin ayant pour label $\w$ est une g\'eod\'esique, on dira que
$\w$ est un {\it mot g\'eod\'esique}.

Pour tout $g\in G$ on pose, si $\ol{h}$ est un sommet de
$\G(G,S)$, $g.\ol{h}=\ol{gh}$, et si $\a$ est l'ar\^ete de label
$s$ entre $\ol{h_1}$ et $\ol{h_2}$, $g.\a$ est l'ar\^ete de label
$s$ entre $g.\ol{h_1}$ et $g.\ol{h_2}$. Celà d\'efinit  une action
à gauche de $G$ sur $\G(G,S)$ ; elle s'effectue par isom\'etrie.\\

 Soit $E$ un espace m\'etrique
g\'eod\'esique. Nous appelons triangle g\'eod\'esique $\lbrack
x,y,z\rbrack$,
 la donn\'ee de 3 points distincts
$x,y,z$ de $E$, et de g\'eod\'esiques les reliant deux \`a deux,
$\lbrack x,y\rbrack, \lbrack y,z\rbrack, \lbrack x,z\rbrack$.
Etant donn\'e un triangle g\'eod\'esique $\lbrack x,y,z\rbrack$ de $E$, il est possible de le plonger isom\'etriquement
dans un triangle $[A,B,C]$ de l'espace euclidien 2-dimensionnel. Appelons $\Psi$ l'isom\'etrie $\Psi :\lbrack
x,y,z\rbrack \longrightarrow \lbrack A,B,C\rbrack$, telle que $\Psi(x)=A,\Psi(y)=B,\Psi(z)=C$.

\begin{figure}[ht]
\center{\includegraphics[scale=0.7]{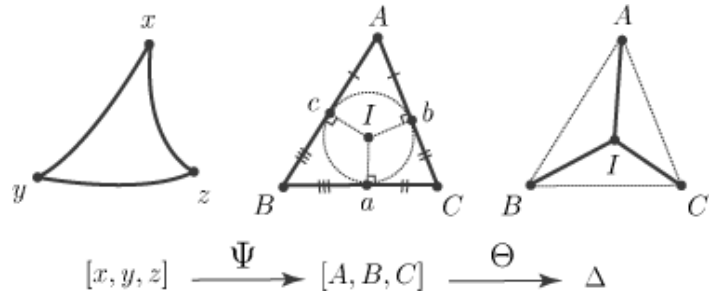}} \caption{Le
trip\^ode $\Delta$ et l'application $\Theta\circ \Psi$.}
\end{figure}

Notons $I$ le centre du cercle inscrit dans $\lbrack A,B,C\rbrack$, et $a,b,c$ ses points de contact (voir la figure
1). Consid\'erons le trip\^ode $\Delta$ constitu\'e des segments $\lbrack A,I\rbrack, \lbrack B,I\rbrack$, $\lbrack
C,I\rbrack$. On construit l'application continue $\Theta$ de $\lbrack A,B,C\rbrack$ dans $\Delta$ de la fa\c con
suivante : la restriction de $\Theta$ sur $\lbrack A,c\rbrack$ est l'unique application affine qui pr\'eserve $A$, et
qui envoie $c$ sur I. On d\'efinit de fa\c con similaire $\Theta$ sur $\lbrack A,b\rbrack$, $\lbrack C,b\rbrack,
\lbrack C,a\rbrack$, $\lbrack B,a\rbrack$, $\lbrack B,c\rbrack$ ; celà d\'efinit $\Theta$ sur $\lbrack x,y,z\rbrack$.
Remarquons que $\Theta$ est bijective sur $A,B,C$, la pr\'e-image de $I$ est $\{a,b,c\}$, et que tout autre point admet
deux pr\'e-images. On construit ainsi l'application $\Theta\circ \Psi$ de $\lbrack x,y,z\rbrack$, qui est unique, \`a
composition par une isom\'etrie du plan pr\`es. Soit $\delta \in \R_+$ ; nous dirons que $E$ est
${\d}$\textit{-hyperbolique}, si pour tout triangle g\'eod\'esique, et pour tout point $m\in \Delta$ (o\`u $\Delta$ est
obtenu par cette construction),
 le diam\`etre de
$(\Theta\circ \Psi)^{-1}(m)$ dans $E$ est major\'e par $\delta$. Nous dirons que $E$ est \textit{hyperbolique} si il
existe $\delta\in \R_+$, tel que $E$ soit $\delta$-hyperbolique.
\\

Un groupe de type fini $G$, muni d'une famille g\'en\'eratrice $S$ sera dit ${\d}$\textit{-hyperboli\-que} si son
graphe de Cayley $\G(G,S)$ est $\d$-hyperbolique. Cette d\'efinition d\'epend du choix de $S$. N\'eanmoins, si $S'$ est
une autre famille g\'en\'eratrice finie pour $G$, et si $\G(G,S)$ est $\d$-hyperbolique, alors $\G(G,S')$ est
$\d'$-hyperbolique pour un certain $\d'$ (en fait à défaut d'être isométrique ils sont {\it quasi-isométriques}, cf.
\cite{cdp}). Nous dirons qu'un groupe est \textit{hyperbolique} si pour une (toute) famille g\'en\'eratrice $S$ finie,
il existe un r\'eel positif $\d$ tel que $\G(G,S)$ soit $\d$-hyperbolique. Ainsi \^etre hyperbolique ne d\'epend pas du
choix d'une famille g\'en\'eratrice finie.

\begin{figure}[!ht]
\centerline{\includegraphics[scale=0.7]{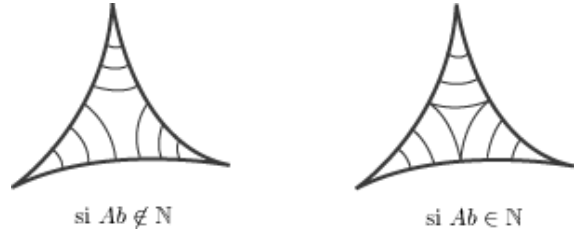}}
\caption{Foliation d'un triangle g\'eod\'esique dans le graphe de
Cayley.}
\end{figure}

Consid\'erons un groupe $\delta$-hyperbolique $G$, et $\G(G,S)$ son graphe de Cayley pour cette famille
g\'en\'eratrice. Soit $\lb x,y,z\rb$ un triangle g\'eod\'esique dans $\G(G,S)$. Soit $\Delta$ le trip\^ode et
l'application $\Theta\circ \Psi$ comme d\'efinis plus haut. Par d\'efinition, si $u,v$ sont deux sommets dans $\lb
x,y,z\rb$, ayant m\^eme image par $\Theta\circ \Psi$, alors ils sont reli\'es dans $\G(G,S)$ par un chemin de longueur
au plus $\d$. La donn\'ee d'un chemin g\'eod\'esique pour tout couple de sommet de $\lb x,y,z\rb$ ayant m\^eme image
par $\Theta\circ
\Psi$, s'appelle une \textit{foliation} de $\lb x,y,z\rb$.\\

Dans la suite on se donne un groupe $G$ hyperbolique, $S$ une famille génératrice finie et $\G=\G (G,S)$ son graphe de
Cayley.

Soit $\g$ un chemin (fini ou infini) dans $\G$ et soient des r\'eels $\l\geq 1,$ $\e\geq 0$. Nous dirons que $\g$ est
une ${(\l,\e)}$\textit{-quasig\'eod\'esique}, si pour tout sous-chemin $\g'$ de $\g$,  si $u$ est le label de $\g'$,
alors,
$$\lgr(u)\leq \l |u|+\e$$
Un chemin sera dit \textit{quasig\'eod\'esique} si c'est une $(\l,\e)$-quasig\'eod\'esique, pour $\l\geq 1,\e\geq 0$.
Nous appelerons quasig\'eod\'esique (resp. bi-)infinie, ou rayon quasigéodésique, un chemin
(resp. bi-)infini qui est une quasig\'eod\'esique.\\

Il est possible de définir le bord d'un espace hyperbolique (\cite{cdp}). Succinctement ici, le bord $\partial \G$ de
$\G$ peut être vu comme l'ensemble des classes d'équivalence de rayons quasigéodésiques par la relation rester à
distance de Hausdorff bornée. On dira qu'une suite de sommets $(s_n)_{n\in\N}$ de $\G$ converge vers $\g\in\partial \G$
si un chemin constitué des géodésiques successives $[s_n,s_{n+1}]$ reste à distance de Hausdorff borné d'un rayon
quasigéodésique dans la classe de $\g$. L'action de $G$ sur $\G$ s'étend naturellement à $\G\cup\partial\G$. On peut
établir que pour tout couple de points $p,q\in\partial \G$ il existe une géodésique bi-infinie (en général non unique)
abusivement notée $[p,q]$ ayant $p,q$ comme points limites dans $\partial \G$.

Nous ferons un usage exclusif et intensif des faits admis qui suivent ; les deux premiers constituent le lemme 10.6.5
et le théorème 3.3.1 de \cite{cdp} :

\begin{lem} Soit $G$ un groupe hyperbolique ; si $g\in G$ est
d'ordre infini alors le chemin bi-infini $\cal{H}=\bigcup_{n\in\Z} h^n.\lb \ol{1},\ol{h}\rb$ est une
quasig\'eod\'esique. Il définit deux points limites $g^-$ et $g^+$ dans $\partial\G$.
\end{lem}

\begin{lem}
Soit $G$ un groupe $\d$-hyperbolique et des réels $\l\geq 1, \e\geq0$. Il existe une constante calculable $k>0$ ne
dépendant que de $\d$, $\l$, $\e$, telle que toute $(\l,\e)$-quasigéodésique bi-infinie ayant pour points limites $h^-$
et $h^+$ dans $\partial\G$ reste à distance de Hausdorff au plus $k$ de toute géodésique $[h^-,h^+]$.
\end{lem}

\begin{figure}[ht]
\center{\includegraphics[scale=0.7]{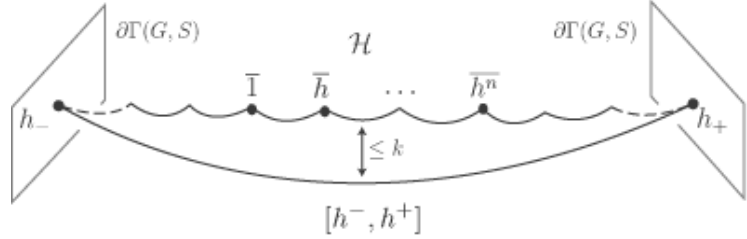}} \caption{La $(\l,\e)$-quasi-g\'eod\'esique $\cal{H}$ ; elle
admet deux points limites $h^-,h^+$ dans $\partial \G(G,S)$, et reste \`a distance de Hausdorff au plus k de toute
g\'eod\'esique $[h^-,h^+]$, o\`u $k$ ne d\'epend que de $\d,\l,\e$.}
\end{figure}

Soit $C\in \N^*$ ; un chemin $\g$ dans $\G$ sera appelé une {\it
$C$-géodésique locale} si tout sous-chemin de $\g$ de longueur au
plus $C$ est une géodésique.
Le lemme suivante est une conséquence immédiate du th\'eor\`eme
III.1.4 de \cite{cdp} (quasi-g\'eod\'esit\'e locale implique
quasi-g\'eod\'esit\'e globale).

\begin{lem}
\label{Creduit} Soit $G$ un groupe muni de la famille
g\'en\'eratrice finie $S$, de fa\c con \`a ce que $(G,d_S)$ soit
$\d$-hyperbolique. Il existe des constantes calculables, $C>0$,
$\l\geq 1$, $\e\geq 0$, telles que tout $C$-géodésique locale est
une $(\l,\e)$-quasigéodésique.
\end{lem}

Enfin, nous ferons un usage intensif des solutions aux problèmes du mot et de la conjugaison dans un groupe
hyperbolique, et en particulier de l'algorithme, dit de Dehn, de réduction monotone d'un mot en mot géodésique (cf.
\cite{cdp}).

\begin{thm}
\label{conj_gro}
 Dans un groupe hyperbolique $G$   :
\begin{itemize}
\item[(i)] On dispose d'un algorithme --nommément l' algorithme de Dehn-- qui donné un mot représentant un \el\ de $G$
retourne un mot géodésique représentant le même \el\ de $G$. En particulier le problème du mot est résoluble dans $G$.
\item[(ii)] Le problème de conjugaison est résoluble dans $G$.
\end{itemize}
\end{thm}

\section{Propri\'et\'e préliminaire}

Nous commençons par d\'emontrer une propri\'et\'e des groupes hyperboliques. Elle nous sera utile tout au long de ce
travail.

\begin{prop}
\label{BS}
Soient $G$ un groupe hyperbolique et $h\in G$ un \el\ d'ordre infini ; il définit $h^-,h^+\in\partial\G$. Le
stabilisateur de $\{ h^-,h^+ \}$ contient $<h>$ comme sous-groupe d'indice fini.
Si de plus il existe $a\in G$ et $p,q\in \Z$, tels que
$$h^q=a h^p a^{-1}\quad dans\; G$$
alors $p=\pm q$ et si $a$ est d'ordre infini $p=q$ ;  $<a>_G$ est
une extension finie d'un sous-groupe de $<h>_G$.
\end{prop}

\noindent \textbf{D\'emonstration} On se place dans le graphe de
Cayley $\G=\G(G,S)$, pour $S$ une famille g\'en\'eratrice finie.
On suppose que pour cette famille $G$ est $\d$-hyperbolique.

Notons $\lb \ol{1},\ol{h}\rb$, un chemin g\'eod\'esique entre
$\ol{1}$ et $\ol{h}$, et consid\'erons le chemin infini $\cal{H}$
d\'efini par $\cal{H}=\bigcup_{n\in\Z} h^n.\lb \ol{1},\ol{h}\rb$.
Remarquons que les sommets $\ol{h^n}$ sont dans $\cal{H}$. Puisque
$h$ est sans torsion, avec le lemme 10.6.5 de \cite{cdp},
$\cal{H}$ est une $(\l,\e)$-quasig\'eod\'esique, avec $\l,\e$ ne
d\'ependant que de $h$. Ainsi $\cal{H}$ admet deux points limites
$h_-$ et $h_+$,
 et si nous notons $\lb h_-,h_+\rb$
une g\'eod\'esique de $h_-$ \`a $h_+$, alors $\lb h_-,h_+\rb$ est
dans le $K$-voisinage de $\cal{H}$, avec $K$ calculable ne
d\'ependant que de $\l,\e$ et $\d$ (th\'eor\`eme 3.3.1 de
\cite{cdp}).

Nous utiliserons le lemme suivant :

\begin{lem}
\label{lemstab} Soient $G$ un groupe $\d$-hyperbolique, et $h\in
G$  tel que $\cal{H}=\bigcup_{n\in\Z} h^n.\lb \ol{1},\ol{h}\rb$
soit une $(\l,\e)$-quasigéodésique de $\G(G)$. Il existe une
constante calculable $K>0$ ne dépendant que de $\d$, $\l$, $\e$
tel que toute classe à gauche de $Stab(h^-,h^+)$ modulo $<h>$
admet un représentant dans $G$ de longueur au plus $K$.
\end{lem}

\noindent{\sl Démonstration du lemme \ref{lemstab}.}
Soit $g\in Stab(h^-,h^+)$ ; $g.\cal H$ est une $(\l,\e)$-quasigéodésique ayant $h^-,h^+$ comme point limite dans
$\partial\G$. Comme nous l'avons vu plus haut $g.\cal H$ et $\cal H$ sont à distance de Hausdorff au plus
$2K(\d,\l,\e)$, avec $K$ calculable, et donc il existe un mot $g'$ de longueur au plus $2K+\lgr(h)$ tel que $g=h^ng'$.
$\ \square$\\

Ce dernier lemme montre que $Stab(h^-,h^+)$ a un nombre fini de
classes à gauche modulo $<h>$, {\it i.e.} $Stab(h^-,h^+)$ contient
$<h>$ comme sous-groupe d'indice fini.

Supposons dans la suite qu'il existe $a\in G$ comme dans l'énoncé. Le calcul montre que pour tout $n\in \Z$,
$$h^{q^n}=a^n h^{p^n} a^{-n}$$
%

Pour $n\in \Z$, consid\'erons $\cal{H}_n=a^n.\cal{H}$, image de $\cal{H}$ sous l'action de $a^n\in G$. Puisque $G$ agit
par isom\'etrie,  pour tout $n$, $\cal{H}_n$ est une $(\l,\e)$-quasig\'eod\'esique. Puisque $h^{q^n}=a^nh^{p^n}a^{-n}$
dans $G$, les quasi-g\'eod\'esiques $\cal{H}$ et $\cal{H}_n$ restent \`a une distance born\'ee (d\'ependant de
$p,q,n,|a|,|h|$). Avec le corollaire 2.1.3 de \cite{cdp}, $\cal{H}_n$ a alors pour points limites $h_-$ et $h_+$, et
avec le th\'eor\`eme 3.3.1 de \cite{cdp},
 pour tout $n$, $\lb h_-,h_+\rb$ est
dans le $K$-voisinage de $\cal{H}_n$.

Consid\'erons un point $\ol{x}$ de $\lb h_-,h_+\rb$. Pour tout $n\in \Z$, il existe un chemin de $\ol{x}$ \`a un point
de $\cal{H}_n$, de longueur au plus $K$ ; et donc pour tout $n$, il existe un entier $i$, et un chemin de $\ol{x}$ \`a
$\ol{a^nh^i}$, de longueur au plus $K+|h|$. Ainsi la boule ferm\'ee centr\'ee en $\ol{x}$ de rayon $K+|h|$, contient
pour tout entier $n$, un sommet de la forme $\ol{a^nh^i}$, pour un certain $i\in \Z$. Puisque le graphe est localement
fini, une boule contient un nombre fini de sommets. Ainsi, il existe $n,m,i,j$, avec $n\not= m$, tels que
$a^nh^i=a^mh^j$ dans $G$. Et donc, il existe $r\in \Z^\ast, s\in \Z$, tels que
$$a^r=h^s$$
Ainsi d'une part $<a>$ est clairement une extension finie de $<h^s>$, et d'autre part $a^r$ et $h$ commutent, et donc,
$$h^{p^r}= a^r h^{p^r} a^{-r}=h^{q^r}$$
ce qui, puisque $h$ est d'ordre infini, n'est possible que si
 $p=\pm q$. De plus
$$h^{qs}=ah^{ps}a^{-1}= a^{rp}=h^{ps}$$
et donc $ps=qs$. Si $a$ est d'ordre infini alors $s\not= 0$ et
nécessairement $p=q$.  $\ \square$\\

\begin{cor}
\label{BS1}
 Soit $G$ un groupe hyperbolique et $h\in G$ un élément d'ordre infini. Le centralisateur et le
normalisateur de $<h>$ contiennent $<h>$ comme sous-groupe
d'indice fini.
\end{cor}

\begin{cor}
Un groupe hyperbolique ne contient pas de groupe de
Baumslag-Solitar comme sous-groupe.
\end{cor}

\begin{cor}
Dans un groupe hyperbolique un \el\ d'ordre infini a une structure de racines presque triviale. En particulier, un
groupe hyperbolique sans torsion a une structure de racines triviale.
\end{cor}

\begin{cor}
Soient $G$ un groupe hyperbolique et $K$ un sous-groupe qui n'est
pas de torsion. Le centralisateur $Z_G(K)$ de $K$ est
virtuellement cyclique. Si $Z_G(K)$ est virtuellement $\Z$ alors
$K$ est aussi virtuellement $\Z$.
\end{cor}

\noindent {\sl Démonstration.} Puisque $K$ n'est pas de torsion il
contient un \el\ $k$ d'ordre infini. $Z_G(K)$ est contenu dans
$Z_G(k)$ qui avec le corollaire \ref{BS1} contient $<k>$ comme
sous-groupe d'indice fini. Ainsi $Z_G(K)$ est virtuellement
cyclique, et infini si et seulement si $Z_G(K)\,\cap <k> \not= \{
1\}$.
 Dans ce dernier cas, il existe $p\not=0$ tel que pour tout $u\in K$, $uk^pu^{-1}=k^p$, et en particulier
 $K\subset Stab(k^-,k^+)$ ; avec la proposition \ref{BS}, $K$
 doit être virtuellement cyclique infini.$\ \square$\\

 Remarquer que $K$ peut être virtuellement $\Z$ sans que
 $Z_G(K)$ le soit (considérer $G=K=\Z_2*\Z_2$).

 \begin{cor}
Soit $G$ un groupe hyperbolique qui n'est pas de torsion et $Z(G)$
son centre. Alors
soit :\\
-- $G$ est élémentaire, $Z(G)$ est virtuellement cyclique,\\
-- $G$ n'est pas élémentaire, $Z(G)$ est fini,\\
-- $G$ n'est pas élémentaire et est sans torsion, $Z(G)$ est
trivial.\\
 \end{cor}

\section{Procédures de réduction cyclique}

On l'a vu : si $h$ est un \el\ sans torsion de $G$ hyperbolique,
alors le chemin bi-infini $\cal{H}=\bigcup_{n\in\Z} h^n.\lb
\ol{1},\ol{h}\rb$ est une $(\l,\e)$-quasig\'eod\'esique. Peut-on,
donn\'e un mot sur les g\'en\'erateurs repr\'esentant $h$,
d\'eterminer algorithmiquement les param\`etres $\l,\e$ de
quasig\'eod\'esit\'e de $\cal{H}$ ?
L'approche que nous
emploierons sera inverse. Donn\'e un  \el\ $h$ sans torsion, nous
le transformerons en un \el\ $h'$ sans torsion tel que
$\bigcup_{n\in\Z} ({h'})^n.\lb \ol{1},\ol{h'}\rb$ soit une
$(\l,\e)$-quasig\'eod\'esique pour $\l,\e$ connus, et tel que $h$
v\'erifie une propri\'et\'e $P$ recherch\'ee, si et seulement si
$h'$ v\'erifie une propri\'et\'e $P'$ analogue. Nous introduirons
pour cela la notion de mot $C$-r\'eduit.\medskip

\noindent \textbf{D\'efinition.} Consid\'erons un groupe $G$
hyperbolique pour une famille g\'en\'eratrice  $S$ donn\'ee. Soit
$C>0$ une constante. Un mot $\w$ sur $S\cup S^{-1}$ sera dit
$C$\textit{-r\'eduit} si $\lgr(\w)>C$ et si $\w$ et tous ses
conjugu\'es cycliques sont des mots g\'eod\'esiques.\smallskip

Donnés les constantes calculables $C$, $\l$, $\e$ provenant du
lemme \ref{Creduit}, si $\w$ est un mot $C$-réduit alors le chemin
infini $\cal{W}=\bigcup_{n\in\Z} \w^n.\lb \ol{1},\ol{\w}\rb$ de
$\G (G,S)$ est une $(\l,\e)$-quasi-g\'eod\'esique ; en
particulier, $\w$ repr\'esente un \el\ sans torsion de $G$.

Consid\'erons un groupe $G$, $\d$-hyperbolique, muni de la famille
g\'en\'eratrice finie $S$. Fixons durant la fin de cette section,
la constante $C$ provenant du lemme \ref{Creduit}. Nous
d\'ecrivons maintenant la proc\'edure permettant donn\'e un mot
sur $S\cup S^{-1}$ repr\'esentant un \el\ $h$ sans torsion de $G$,
de lui associer un mot $C$-r\'eduit repr\'esentant un conjugu\'e
d'une puissance de $h$.\\

\noindent \textbf{Proc\'edure de $C$-r\'eduction :}
 Donné un mot $\w$ sur $S\cup S^{-1}$ repr\'esentant un \el\ non trivial $h\in
 G$,
 la proc\'edure de $C$-r\'eduction retourne, dès-lors que $h$ est d'ordre infini, un 3-uplet
 $(h_1,u,n)\in S^*\times S^*\times \N$ tel que
$h_1$ est un mot $C$-r\'eduit et
$h_1=uh^nu^{-1}$ dans $G$.\smallskip\\
Initialement on applique l'étape 1 au 3-uplet $(\w,1,1)\in
(S^*,S^*,\N)$.\\
\noindent\textbf{Etape 1 :} Donné $(\w,u,n)$ on applique d'abord l'algorithme de Dehn pour transformer
 $\w$ en un mot g\'eod\'esique.
 Un conjugu\'e cyclique de $\w$ peut ne
pas être g\'eod\'esique ; par exemple $\w\equiv \w_a\w_b$ et il existe un mot g\'eod\'esique $\w_1$, tel que
$$\w_b\,\w_a=\w_1\quad \mathrm{dans}\; G$$
et
$$\lgr(\w_1)<\lgr(\w_b\,\w_a)=\lgr(\w_b)+\lgr(\w_a)=\lgr(\w)$$
On applique l'algorithme de Dehn  à tous les conjugués cycliques de $\w$ pour en d\'ecider et le cas \'ech\'eant
d\'eterminer $\w_1$ comme ci-dessus ; il repr\'esente dans $G$ le conjugu\'e de $\w$ par $\w_a^{-1}$. On applique la
même procédure à $\w_1$, pour obtenir \'eventuellement $\w_2$, et ainsi de suite. Puisque la longueur des mots diminue
strictement, la proc\'edure s'arr\^ete, et on obtient un mot $\w_0$ non vide, et un mot $v$, tels que tous les
conjugu\'es cycliques de $\w_0$ soient g\'eod\'esiques, et tels que $\w_0=v\w v^{-1}$ dans $G$. Si $\w_0$ est de
longueur au moins $C$, la proc\'edure s'arr\^ete et retourne $(\w_0,vu,n)$,
sinon on applique l'\'etape 2 \`a $(w_0,vu,n)$.\smallskip\\
\noindent \textbf{Etape 2 :} Donné $(\w,u,n_0)$ où $\w$ est un mot de longueur inf\'erieure \`a $C$,  on consid\`ere
les puissances successives de $\w$ et on applique l'algorithme de Dehn pour les transformer en mots géodésiques jusqu'à
trouver un entier $n$ tel que  dans $G$, $|w^n|\geq C$. Si $h$ est d'ordre infini la procédure s'arrête (car la boule
de rayon $C$ contient un nombre fini d'éléments) et on applique alors l'\'etape 1 \`a $(\w^n,u,nn_0)$ ; sinon on finit
par trouver $n$ tel que $\w^n=1$ dans $G$.
\\

Dans le cas o\`u nous traiterons avec deux éléments $h,h'$, nous
aurons besoin de raffiner la proc\'edure de
$C$-r\'eduction.\smallskip

\noindent $\textbf{Proc\'edure de\;} {C}$\textbf{-r\'eductions jumel\'ees :} Donnés deux mots $h,h'$ repr\'e\-sentant
des \el s d'ordre infini de $G$, la proc\'edure retourne deux mots $C$-réduits $h_1,h_1'$, $u,u'\in G$ et un entier
$n\in\N_*$ tels que $h_1$ et $h_1'$ soient conjugu\'es respectivement  \`a $h^n$ par $u$ et à ${h'}^n$ par $u'$. Pour
cel\`a on adapte la procédure de $C$-réduction de façon à l'appliquer à deux 3-uplets $(h,u,n)$ et $(h',u',n')$ de
$S^*\times S^*\times \N$
simultan\'ement,  de la fa\c con suivante :\smallskip\\
On se donne initialement les deux 3-uplets $(h,1,1)$ et $(h',1,1)$
de $S^*\times S^*\times \N$, auxquels on applique l'étape 1.\\
 \noindent \textbf{Etape 1 :} Il s'agit de l'étape 1 de l'algorithme de $C$-réduction
 appliqué en parallèle à deux 3-uplets $(h,u,n)$ et $(h',u',n')$. A la fin de cette étape, on
retourne un résultat seulement si l'on obtient des
 représentants $C$-réduits des classes de conjugaison et de $h$ et de $h'$ et sinon on passe à l'étape 2.
 \smallskip\\
\noindent\textbf{Etape 2 :} C'est l'étape 2 de l'algorithme de $C$-réduction appliqué simultanément aux deux 3-uplets
$(h,u,n)$ et $(h',u',n')$, et ce jusqu'à d\'eterminer  un $N\in \N_\ast$, tel que à la fois $|h^N|\geq C$ et
$|{h'}^N|\geq C$ ; on passe alors à l'étape 1.\smallskip\\
\indent Le m\^eme argument que
 dans de le cas de la $C$-r\'eduction montre que la proc\'edure s'arr\^ete.\\

\section{Les algorithmes}

Maintenant que nous avons d\'efini les notions d'\el\ $C$-r\'eduit et de proc\'edure de $C$-r\'eduction nous allons
apporter des solutions à divers probl\`emes algorithmiques  dans les groupes hyperboliques. Nous utiliserons de fa\c
con syst\'ema\-tique les solutions aux probl\`emes du mot et de la conjugaison.
Nous consid\'erons dans la suite un groupe $G$ muni d'une famille g\'en\'eratrice finie $S$ avec $(G,d_S)$
$\d$-hyperbolique et nous fixons les constantes $C,\l,\e$, provenant du lemme \ref{Creduit}.

\subsection{Déterminer la torsion}
Comme conséquence directe du paragraphe précédent, on obtient :

\begin{thm}[\textit{D\'eterminer la torsion}]
\label{torsion} Un groupe hyperbolique $G$ contient un nombre fini
de classes de conjugaison d'éléments de torsion.

Donné un mot sur les générateurs on peut décider algorithmiquement
s'il représente ou non un élément de torsion de $G$.

On peut algorithmiquement se donner un représentant de chaque
classe de conjugaison d'élément de torsion.
\end{thm}

\noindent \textsl{D\'emonstration.} Pour décider si un mot
représente un élément de torsion il suffit de lui appliquer la
procédure de $C$-réduction. On finira par trouver soit un conjugué
$C$-réduit (et donc d'ordre infini) soit son ordre (fini).

Tout élément de torsion a un conjugué dans la boule centrée en 1
de rayon $C$. Cette dernière ne contient qu'un nombre fini
d'éléments. Ainsi d'une part $G$ contient un nombre fini de
classes de conjugaison d'éléments de torsion, et d'autre part, en
décidant pour chaque élément de $B(1,C)$ s'il est d'ordre infini
on détermine un représentant
pour chaque classe de conjugaison d'éléments de torsion.$\ \square$\\

\subsection{Sous-groupes cycliques et problème du mot généralisé}

\begin{lem}
\label{V1} Soit $G$ un groupe $\d$-hyperbolique, muni de la
famille g\'en\'eratrice finie, $S$. Soient $\d\geq 1$, $\e\geq 0$,
et $C$ comme dans le lemme \ref{Creduit}. Soient $h$ et $\w$ deux
mots sur $S$, $C$-r\'eduits. Si $\w$ est conjugu\'e dans $G$ \`a
$h^n$, alors il existe une constante calculable $L=L(\d,\e,\l,|\w
|,|h |)$, tel que $n\leq L$.
\end{lem}

\noindent \textsl{D\'emonstration} On se place dans le graphe de
Cayley. On note $\lb \ol{1},\ol{h}\rb$ le chemin g\'eod\'esique de
label $h$, et $\mathcal{H}$ le chemin bi-infini :
$$\mathcal{H}=\bigcup_{n\in\Z}h^n.\lb\ol{1},\ol{h}\rb$$
Puisque $h$ est $C$-r\'eduit, avec le lemme \ref{Creduit},
$\mathcal{H}$ est une $(\l,\e)$-quasig\'eod\'esique de $\G (G,S)$.
Avec le th\'eor\`eme III.3.1 de \cite{cdp}, (th\'eor\`eme de
stabilit\'e des quasi-g\'eod\'esiques de longueur infinie),
$\mathcal{H}$
 a
exactement deux points d'accumulation, $h_-$ et $h_+$ dans
$\partial \G (G,S)$, et il existe une constante calculable $k$, ne
d\'ependant que de $\d,\l,\e$, telle que $H$ est dans un
$k$-voisinage de toute g\'eod\'esique joignant $h_-$ \`a $h_+$,
et inversement, toute g\'eod\'esique joignant $h_-$ \`a $h_+$ est
dans le $k$-voisinage de $H$.

Notons $\lb \ol{1},\ol{\w}\rb$, le chemin g\'eod\'esique
de label $\w$. Supposons que $h^n=u \w u^{-1}$ dans $G$. Puisque $\w$ est
$C$-r\'eduit, et que $G$ agit par isom\'etrie
sur $\G (G,S)$, le chemin bi-infini,
$$\mathcal{W}=u.(\bigcup_{n\in\Z}\w^n.\lb\ol{1},\ol{\w}\rb )$$
est une $(\l,\e)$-quasi-g\'eod\'esique, et a donc deux points d'accumulation dans $\partial \G (G,S)$. Maintenant,
puisque $h^n=u \w u^{-1}$, $\mathcal{H}$ et $\mathcal{W}$ restent \`a distance de Hausdorff born\'ee. Avec le
corollaire II.1.3 de \cite{cdp}, $\mathcal{W}$ a donc $h_-$ et $h_+$ comme points d'accumulation.

Consid\'erons maintenant, $[h_-,h_+]$, une g\'eod\'esique reliant
$h_-$ et $h_+$. Avec ce qui pr\'ec\`ede, $\mathcal{H}$ est
dans le $k$-voisinage
de $[h_-,h_+]$ qui est elle-m\^eme dans le $k$-voisinage
\begin{figure}[ht]
\center{\includegraphics[scale=0.7]{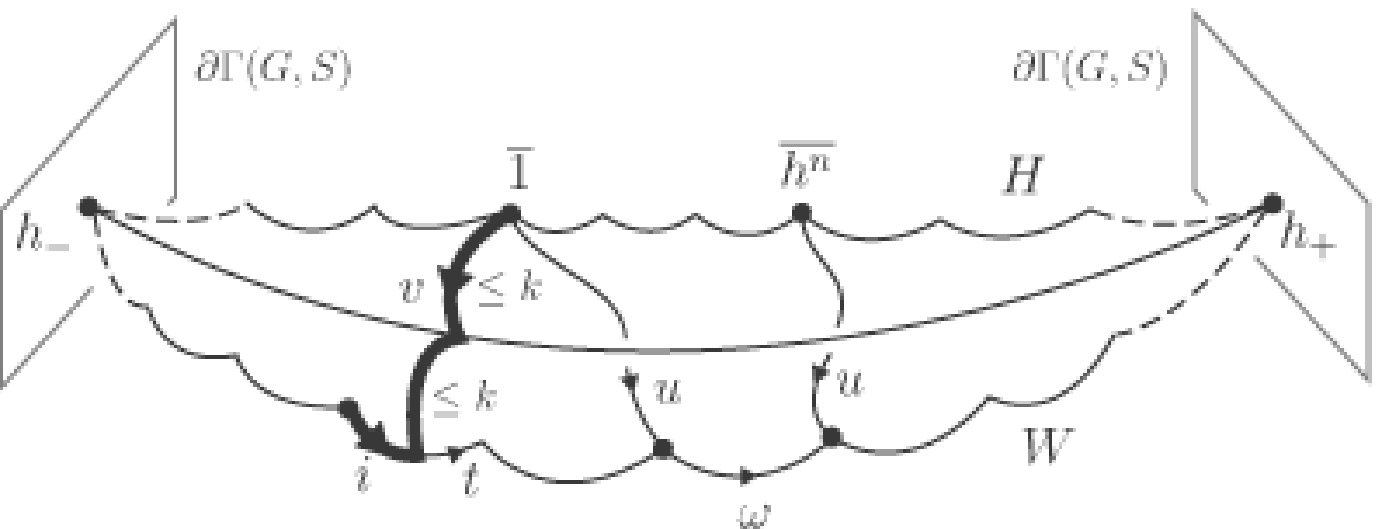}} \caption{}
\end{figure}
 de $W$.
Ainsi il existe un chemin reliant $\ov{1}$ \`a un point de
$\mathcal{W}$, de longueur inf\'erieure \`a $2k$ (cf. figure 4).
 Notons $v$ le label de ce
chemin. Alors, il existe $p\in \Z$, tel que dans $G$,
$$v=u \w^p \w_i$$
o\`u $\w_i$ est un sous-mot initial de $\w$, \emph{i.e.} $\w\equiv \w_i.\w_t$.
Ainsi, $u=v\w_i^{-1}\w^{-p}$, et puisque $h^n=u\w u^{-1}$,
\begin{align*}
h^n &= v\w_i^{-1}\w^p \w \w^{-p} \w_i v^{-1}\\
&= (v\w_i^{-1}) \w (v\w_i^{-1})^{-1}\\
&= v' \w {v'}^{-1}
\end{align*}

On obtient donc un \el\ $v'=v\w_i^{-1}$ de longueur major\'ee par $2k + |w|$,
qui conjugue $\w$ en $h^n$. Maintenant, puisque $\mathcal{H}$ est une
$(\l,\e)$-quasi-g\'eod\'esique,
\begin{align*}
\lgr(h^n)&\leq \l |h^n| + \e \\
n \lgr(h) &\leq \l (4k+3|\w |)+\e\\
n &\leq \frac{1}{\lgr(h)}(\l (4k+3|\w |)+\e)
\end{align*}
et ainsi, $n$ est major\'e par une constante $L$, ne d\'ependant que de $\d,\l,\e,|\w|,|h|$. $\ \square$

\begin{thm}[\textit{Problème du mot généralisé à conjugaison près de $\Z$ dans
$G$}] \label{conj_per_gro}

Soit $G$ un groupe hyperbolique, et $h\in G$ un \el\ d'ordre
infini. Donn\'e un \el\ $\w\in G$, on peut d\'ecider si $\w$ est
conjugu\'e dans $G$ \`a un \el\ de $<h>$. De plus, si c'est le
cas,  $\w$ est conjugu\'e \`a au plus deux \el s (resp. un \el\ si
$G$ est sans torsion) de $<h>$.
\end{thm}

\noindent \textsl{D\'emonstration.} On considère une famille
génératrice $S$ pour laquelle $(G,d_S)$ soit $\d$-hyperbolique. On
suppose $h$ et $\w$ donn\'es par des mots sur $S$. On fixe $\l\geq
1$ et $\e\geq 0$, et $C$ donn\'es par le lemme \ref{Creduit}. Si
les mots $h$ et $\w$ sont $C$-r\'eduits, alors le lemme \ref{V1}
nous donne une borne calculable $L$, ne d\'ependant que de $\d ,
\l ,\e ,|\w |,|h|$, telle que si $\w$ est conjugu\'e \`a $h^n$,
pour $n\in \Z$, alors $|n|\leq L$. On utilise ensuite l'algorithme
de la conjugaison pour d\'ecider si $\w$ est conjugu\'e \`a $h^p$,
pour tout entier $p$ v\'erifiant $-L\leq p\leq L$. On d\'etermine
ainsi tous les \el s $p\in \Z$, tels que $\w$ soit conjugu\'e \`a
$h^p$ (avec la proposition \ref{BS}, il y en a au plus deux, un
seul si $G$ n'a pas de torsion).

Si $h$ ou $\w$ ne sont pas $C$-r\'eduits, alors on leur
applique la proc\'edure de $C$-r\'eduction jumel\'ee, d'une part,
et d'autre part, on utilise l'algorithme du mot pour
tester si une puissance de $\w$ repr\'esente l'identit\'e de $G$.
Si $\w$ est de torsion, cet algorithme finit par le d\'eterminer,
et $\w$ ne peut d\'efinitivement pas \^etre conjugu\'e \`a une
puissance non triviale de $h$. Sinon, on finit par obtenir
les mots $h_1,\w_1$ $C$-r\'eduits, conjugu\'es respectifs
de $h^N$ et $\w^N$ pour un certain $N$.

Si $\w$ est conjugu\'e \`a $h^n$, alors $\w_1$ est conjugu\'e \`a $h_1^n$. En utilisant le lemme \ref{V1}, on
d\'etermine au plus deux \el s $n\in \Z$ tel que $\w_1$ soit conjugu\'e \`a $h_1^n$. Il suffit pour conclure d'utiliser
l'algorithme de la conjugaison, pour d\'ecider si $\w$ est conjugu\'e \`a $h^n$ pour un tel $n$.  $\ \square$

\begin{thm}[\textit{Probl\`eme du mot g\'en\'eralis\'e de}\ $\Z$\
\textit{dans}\ $G$]
 \label{gwp_gro} Soient $G$ un
\linebreak groupe hyperbolique, et $h\in G$ un \el\ d'ordre
infini. Donn\'e $\w\in G$ on peut d\'ecider si $\w\in <h>$.
\end{thm}

\noindent \textsl{D\'emonstration.} Appliquer l'algorithme du mot, pour d\'ecider si $\w$ est \'egal \`a un des \el s
de $<h>$ qui lui sont conjugu\'es, fournis par l'algorithme du théorème \ref{conj_per_gro}.  $\ \square$

\begin{thm}[{\it Problème du mot généralisé d'un sous-groupe virtuellement $\Z$}]
Soit $G$ un groupe hyperbolique et $H$ un sous-groupe
virtuellement $\Z$ de $G$. Le problème du mot généralisé de $H$
dans $G$ est résoluble.
\end{thm}

\noindent \textsl{Démonstration.} Supposons que $H=H_0\cup
h_1H_0\cup\cdots \cup h_nH_0$ avec $H_0$ cyclique infini ; pour
décider si $\w\in H$ il suffit de décider avec l'algorithme du
théorème \ref{gwp_gro} si $\w$, ou $h_1^{-1}\w$, ..., ou
$h_n^{-1}\w$ est dans
$H_0$.$\ \square$\\

\subsection{Malnormalité de sous-groupes cycliques infinis}

\begin{lem}
\label{V2} Soit $G$ un groupe $\d$-hyperbolique, muni d'une
famille g\'en\'eratrice finie $S$. Soient $\l\geq 1$, $\e\geq 0$,
et $C$ la constante dont l'existence provient du lemme
\ref{Creduit}. Soient $h_1$ et $h_2$ des mots sur $S$,
$C$-r\'eduits, et $u$ et $v$ des \el s de $G$.

Supposons qu'il existe $n_1,n_2\in \Z$ tels que
$u=h_1^{n_1}vh_2^{n_2}$ dans $G$. Alors il existe une constante calculable
$K>0$, ne d\'ependant que de $\d,\l,\e,|u|,|v|,|h_1|,|h_2|$ et du
cardinal de $S$, et deux entiers $m_1,m_2$ avec
$|m_1|<K$ ,$|m_2|<K$, tels que $u=h_1^{m_1}vh_2^{m_2}$ dans $G$.

\end{lem}

\noindent \textsl{D\'emonstration.} On se place dans le graphe de
Cayley $\G (G,S)$, et on suppose que $u=h_1^{n_1}vh_2^{n_2}$.
Notons $\lb \ol{1}, \ol{h_1}\rb$, et $\lb\ol{1},\ol{h_2}\rb$ les
chemins g\'eod\'esiques de labels respectifs $h_1,h_2$. On note
$$\mathcal{H}_1=\lb\ol{1},\ol{h_1}\rb\cup h_1.\lb\ol{1},\ol{h_1}\rb\cup\cdots
\cup h_1^{n_1-1}.\lb\ol{1},\ol{h_1}\rb$$
$$\mathcal{H}_2=\lb\ol{1},\ol{h_2^{-1}}\rb\cup h_2^{-1}.\lb\ol{1},\ol{h_2^{-1}}\rb
\cup \cdots \cup h_2^{1-n_2}.\lb\ol{1},\ol{h_2^{-1}}\rb$$

Puisque $h_1$ et $h_2$ sont des mots $C$-r\'eduits, avec le lemme
\ref{Creduit}, ce sont des chemins $(\l,\e)$-quasig\'eod\'esiques
de longueur finie. Puisque $G$ agit par isom\'etrie sur $\G
(G,S)$, le chemin $u.\mathcal{H}_1$ est aussi une
$(\l,\e)$-quasi-g\'eod\'esique.

\begin{figure}[ht]
\center{\includegraphics[scale=0.7]{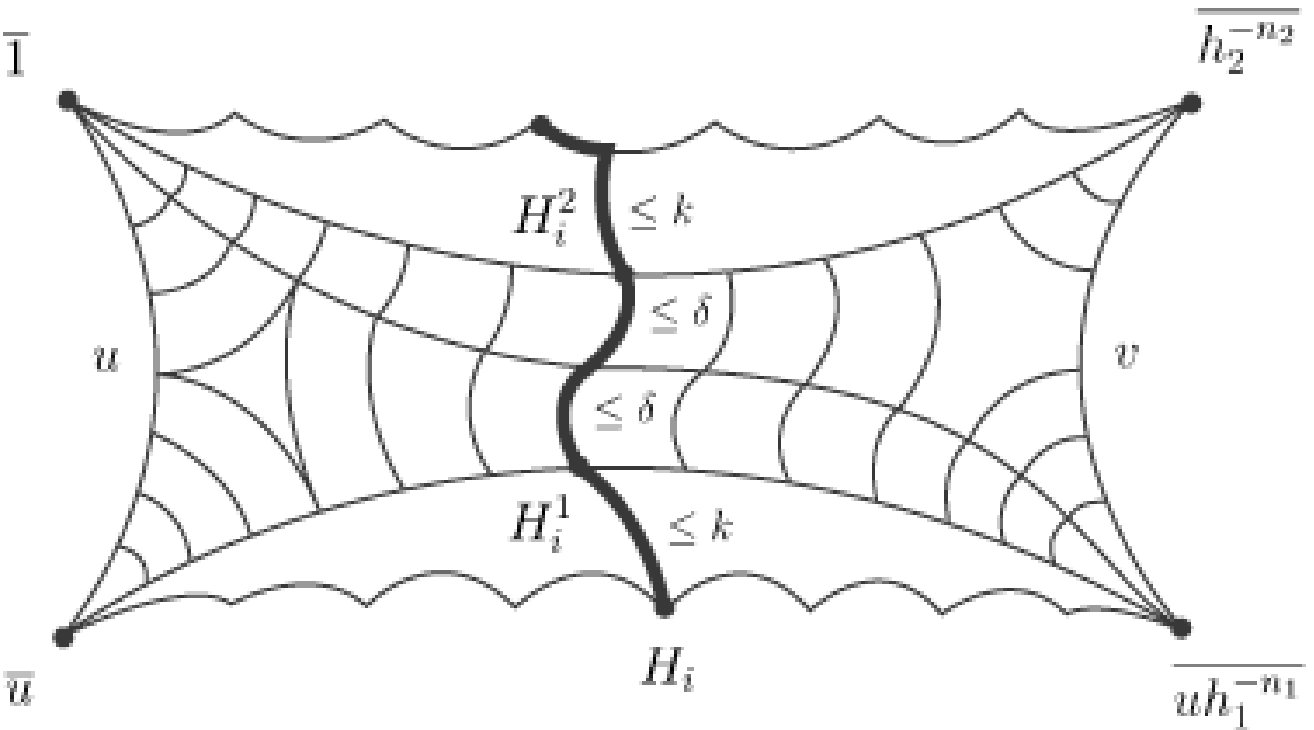}} \caption{}
\end{figure}

Avec le th\'eor\`eme de stabilit\'e des quasi-g\'eod\'esiques de
longueur finie (th\'eor\`eme III.1.3 de \cite{cdp}, il existe une
constante calculable $k$ ne d\'ependant que de $\d,\l,\e$, telle que
$u.\mathcal{H}_1$, est dans le $k$-voisinage de toute
g\'eod\'esique joignant $\ol{u}$ \`a $\overline{uh_1^{n_1}}$.

Consid\'erons un quadrilat\`ere  g\'eod\'esique $Q$, de sommets
$\ol{1},\ol{u},\overline{uh_1^{n_1}},\overline{h_2^{-n_2}}$.
Soit $0\leq n\leq n_1$, et notons $H_n$ le point
$\overline{uh_1^n}$. Il existe un chemin de $H_n$ \`a
la g\'eod\'esique $[\ol{u},\overline{uh_1^{n_1}}]$. Appelons
$H^1_n$ son extr\'emit\'e. En consid\'erant une g\'eod\'esique
de $\ol{1}$ \`a $\overline{uh_1^{n_1}}$, puis une foliation des
deux triangles g\'eod\'esiques obtenus \`a partir de $Q$ et de
cette
g\'eod\'esique, on peut construire un chemin reliant
$H_n^1$ \`a $[\ol{1},\overline{h_2^{-n_2}}]$, de longueur
inf\'erieure \`a $k_1=3\d + \max\{|u|,|v|\}$. Appelons
$H_n^2$ son extr\'emit\'e. Puisque $\mathcal{H}_2$ est
une $(\l,\e)$-quasi-g\'eod\'esique, avec le th\'eor\`eme
de stabilit\'e des quasi-g\'eod\'esiques de longueur
finie,
il existe un chemin de $H_n^2$ \`a $\mathcal{H}_2$, de longueur
inf\'erieure \`a $k$. Ainsi, on peut trouver
$0\leq m\leq n_2$, et un chemin de $H_n^2$ \`a $\ol{h_2^{-m}}$, de
longueur inf\'erieure \`a $k+|h_2|$. Et finalement, pour tout $n,
0\leq n\leq n_1$, On peut trouver un chemin de $H_n$ \`a
$\ol{h_2^{-m}}$ de longueur major\'ee par $K_1=2k + 3\d
+\max\{|u|,|v|\}+|h_2|$ (cf. figure 5).

Il existe au plus $\mathrm{card}(S)^{K_1}$ \el s de $G$ pouvant
s'\'ecrire en un mot de longueur au plus $K$. Supposons que $n_1>
\mathrm{card}(S)^{K_1}$. Alors il existe $c\in G$, $0<r<n_1$, et
$s\in \Z$ tels que $h_1^r= c h_2^s c^{-1}$ dans $G$. Ainsi
$u=h_1^{n_1-r} v h_2^{n_2+s}$, avec $n_1-r< n_1$ (cf. figure 6).
 Ainsi on
peut trouver $m_1,m_2$, avec $m_1\leq \mathrm{card}(S)^{K_1}$, tels que
$u=h_1^{m_1}vh_2^{m_2}$.

Il ne reste plus qu'\`a majorer $m_2$ (on ne peut pas utiliser un
argument par sym\'etrie, car lorsque $n_1$ diminue $n_2$ peut
augmenter, et r\'eciproquement). Le chemin
$[\ol{1},\ol{h_2^{-1}}]\cup\cdots \cup
h_2^{1-n_2}[\ol{1},\ol{h_2^{-1}}]$ est une
$(\l,\e)$-quasi-g\'eod\'esique, ainsi,
\begin{align*}
  \lgr(h_2^{m_2}) &\leq \l |h_2^{-m_2}|+\e \\
  m_2 |h_2| &\leq \l (|u|+|v|+ m_1|h_1|)+\e \\
  m_2 &\leq \frac{1}{|h_2|} (\l (|u|+|v|+ m_1|h_1|)+\e)
\end{align*}

\begin{figure}[ht]
\center{\includegraphics[scale=0.7]{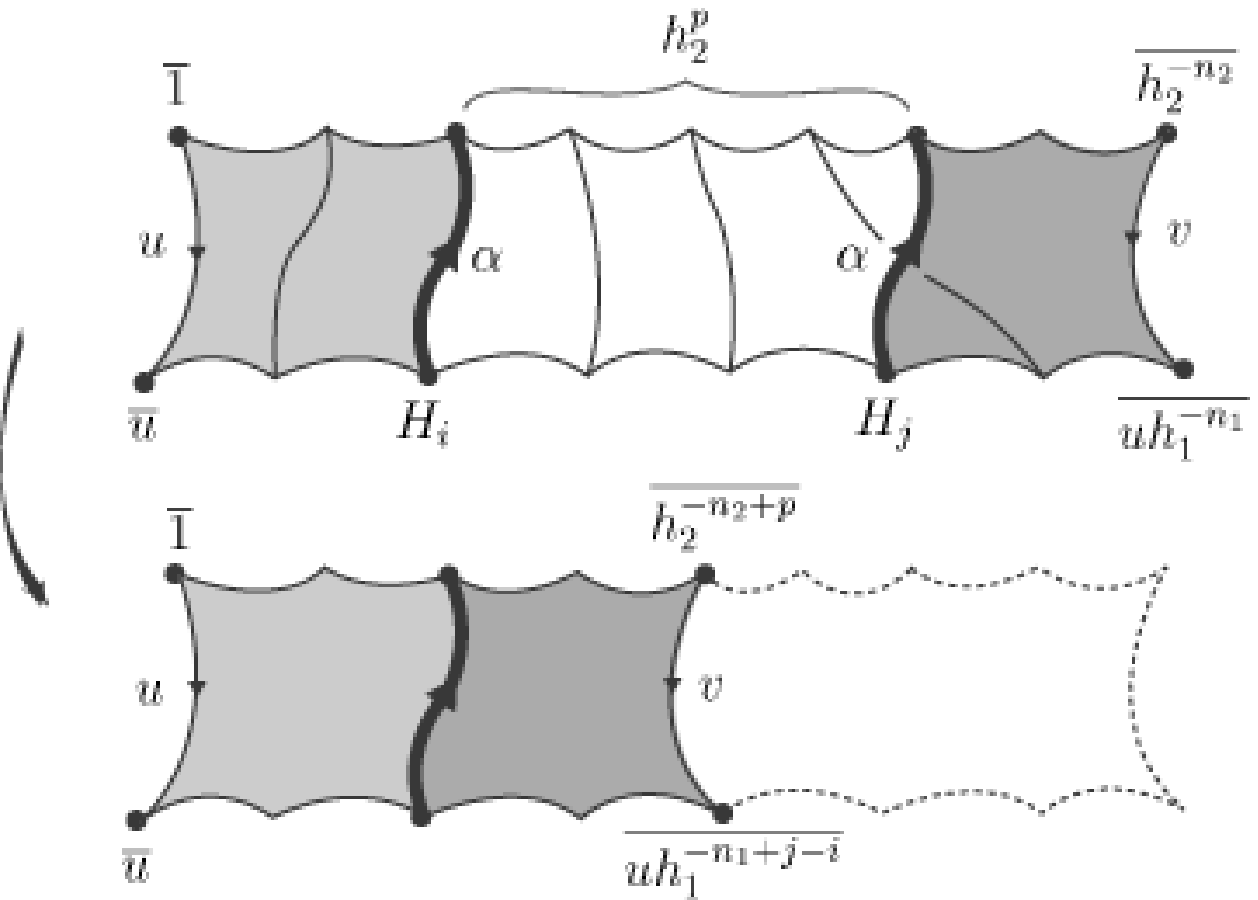}} \caption{}
\end{figure}

Ainsi $m_2$ est major\'e par une constante $K_2$ ne d\'ependant que de $\d,\l,\e,|u|,|v|$,$|h_1|$,\linebreak et
$|h_2|$. Maintenant on pose $K=\max\{\mathrm{card}(S)^{K_1}, K_2\}$. $\ \square$

\begin{thm}
\label{2coset_gro} Soit $G$ un groupe hyperboli\-que et soient
$h_1$ et $h_2$ des \el s de $G$ d'ordre infini. Donn\'es $u,v\in
G$, on peut d\'ecider si il existe $n_1,n_2\in \Z$ tels que
$u=h_1^{n_1}vh_2^{n_2}$, et le cas \'ech\'eant donner un tel
couple $(n_1,n_2)$
\end{thm}

\noindent \textsl{D\'emonstration.} On se donne une famille
génératrice $S$ telle que $(G,d_S)$ soit $\d$-hyperbolique ainsi
que $\l\geq 1$, $\e\geq 0$, et $C$ comme dans le lemme
\ref{Creduit}. On se donne $h_1, h_2, u, v$ par des mots sur $S$.
Si $h_1$ et $h_2$ sont $C$-r\'eduits, le lemme \ref{V2} nous donne
une constante calculable $K$, et on peut alors utiliser
l'algorithme du mot pour d\'eterminer, si il existe, deux entiers
$m_1,m_2$, avec $|m_1|<K$, $|m_2|<K$, tels que
$u=h_1^{m_1}vh_2^{m_2}$.

Si $h_1$ ou $h_2$ n'est pas $C$-r\'eduit,
on applique le proc\'ed\'e de $C$-r\'eduction jumel\'ee \`a $h_1$ et $h_2$.
On obtient des mots $H_1,H_2$,
$C$-r\'eduits, les mots $a_1,a_2$, et $N\in \N$,
tels que $H_1=a_1h_1^Na_1^{-1}$ et
$H_2=a_2h_2^Na_2^{-1}$ dans $G$.
Alors, $n_1= q_1 N + r_1$ et $n_2=q_2 N +r_2$, pour
$q_1,q_2\in \Z$, et $r_1,r_2\in \{0,\ldots ,N-1\}$, et
\begin{align*}
u &= h_1^{n_1}v h_2^{n_2}\\
 &= h_1^{q_1N+r_1}vh_2^{q_2N+r_2}\\
\Leftrightarrow\quad uh_2^{-r_2} &= (h_1^N)^{q_1} h_1^{r_1}v
 (h_2^N)^{q_2}\\
\Leftrightarrow \quad a_1uh_2^{-r_2}a_2^{-1}&=
 a_1(h_1^N)^{q_1}a_1^{-1} a_1h_1^{r_1}va_2^{-1}
 a_2(h_2^N)^{q_2}a_2^{-1}\\
\Leftrightarrow \quad a_1uh_2^{-r_2}a_2^{-1}&=
 {H_1}^{q_1} (a_1h_1^{r_1}va_2^{-1})
 {H_2}^{q_2}
\end{align*}
Pour $0\leq i< N$, on pose
$$U_i=a_1uh_2^{-i}a_2^{-1}$$
$$V_i=a_1h_1^i v a_2^{-1}$$

Ainsi il existe $n_1,n_2\in \Z_\ast$, tels que $u=h_1^{n_1}vh_2^{n_2}$, si et seulement si il existe $q_1,q_2\in \Z$,
et $i,j\in\{0,\ldots ,N-1\}$ tels que $U_i={h_1'}^{q_1} V_j {h_2'}^{q_2}$, ce dont on peut d\'ecider puisque
$\{0,\ldots ,N-1\}$ est fini, et que $H_1,H_2$ sont $C$-r\'eduits. $\ \square$

\begin{lem}
\label{V3} Soit $G$ un groupe $\d$-hyperbolique muni d'une famille
g\'en\'eratrice $S$. Soient $\d,\e$ et $C$ comme dans le lemme
\ref{Creduit}, et soient $h_1,h_2$ des mots sur $S$,
$C$-r\'eduits. Si un \el\ non trivial de $<h_1>$ est conjugu\'e
dans $G$ \`a un \el\ de $<h_2>$, alors il existe une constante
$L=L(\d,\e,\l,|h_1|,|h_2|)$, et $n_1,n_2\in {\Z}_\ast$ tels que
$|n_1|\leq L$, $|n_2|\leq L$, et $h_1^{n_1}$ est conjugu\'e \`a
$h_2^{n_2}$.
\end{lem}

\noindent \textsl{D\'emonstration.} Supposons qu'il existe $u\in
G$, qui conjugue $h_1^{p_1}$ en $h_2^{p_2}$.
Notons $[\ol{1},\ol{h_1}]$ et $[\ol{1},\ol{h_2}]$ les chemins
g\'eod\'esiques de labels respectifs $h_1,h_2$. Avec le lemme
\ref{Creduit}, puisque $h_1,h_2$ sont $C$-r\'eduits, les chemins
bi-infinis de $\G (G,S)$,
$\mathcal{H}_1=\bigcup_{n\in\Z}h_1^n.[\ov{1},\ol{h_1}]$, et
$\mathcal{H}_2=\bigcup_{n\in\Z}h_2^n.[\ol{1},\ol{h_2}]$ sont des
$(\l,\e)$-quasi-g\'eod\'esiques. Puisque $u$ conjugue $h_1^{p_1}$
en $h_2^{p_2}$, $\mathcal{H}_1$ et $\mathcal{H}_2$
 sont \`a
distance de Hausdorff born\'ee, et donc avec  le th\'eor\`eme III.3.1, et le corollaire II.1.3 de \cite{cdp}, elles ont
m\^eme points d'accumulation dans le bord du graphe de Cayley, $h_-$ et $h_+$. Avec le th\'eor\`eme III.3.1 de
\cite{cdp} elles sont dans un $2k$-voisinage l'une de l'autre, o\`u $k$ ne d\'epend que de $\d,\l,\e$. De la m\^eme
fa\c con que dans la d\'emonstration du lemme \ref{V1}, on peut alors trouver un \el\ $v\in G$ de longueur inf\'erieure
\`a $2k+|h_1|+|h_2|$, qui conjugue $h_1^{p_1}$ en $h_2^{p_2}$. En appliquant le lemme \ref{V2} il existe donc une
constante $L=L(\d,\e,\l,|h_1|,|h_2|)$, et $n_1,n_2\in \Z$ tels que $|n_1|\leq L$, $|n_2|\leq L$, et $h_1^{n_1}$ est
conjugu\'e \`a $h_2^{n_2}$. $\ \square$

\begin{thm}[{\it Malnormalité d'un sous-groupe cyclique infini}]
\label{lim_gro_per}
Soient $G$ un groupe hyperbolique et $h$ un \el\ d'ordre infini de $G$. Il existe un algorithme permettant de d\'ecider
s'il existe un \el\ de $G\, -<h>$ qui conjugue deux \el s non triviaux de $<h>$. En particulier on peut décider si
$<h>$ est malnormal dans $G$.
\end{thm}

\noindent \textsl{D\'emonstration.} Commen\c cons par remarquer,
bien que cela ne soit pas n\'ecessaire,
 qu'avec la proposition \ref{BS},
si $a$ conjugue $h^p$ en $h^q$, alors n\'ecessairement $p=\pm q$
et $a^r=h^s$.

Supposons tout d'abord que $h$ soit $C$-r\'eduit. Avec le lemme
\ref{V3}, on peut prendre $0< p\leq L$ o\`u $L$ est un entier
calculable ne d\'ependant que de $\d,\l,\e,|h|$. Reprenons la
d\'emonstration du lemme \ref{V1}. Si $u$ conjugue $h^p$ en $h^q$,
alors il existe $v'$ qui conjugue $h^p$ en $h^q$ v\'erifiant en
outre $|v'|\leq 2k + |h|$ (o\`u $k$ ne d\'epend que de
$\d,\l,\e$). De plus,  $u=v' h^t$, et ainsi, $u\in <h>\ \iff\
v'\in <h>$. Ainsi il suffit d'utiliser l'algorithme du mot pour
d\'ecider si $h^p=uh^qu^{-1}$, pour $0< p\leq L$, $p=\pm q$, et
$|u|\leq 2k +|h|$, et l'algorithme du théorème \ref{gwp_gro} pour
d\'ecider si $u\in <h>$, pour conclure.

Si $h$ n'est pas $C$-r\'eduit, on le r\'eduit en un \el\ $C$-r\'eduit $h'=ah^na^{-1}$. Il est ais\'e de v\'erifier que
s'il existe un \el\ hors de $<h>$, qui conjugue deux \el s non triviaux de $<h>$, alors il existe un \el\ hors de
\mbox{$<aha^{-1}>\supset <h'>$} (de longueur born\'ee), qui conjugue deux \el s non triviaux de $<h'>$.
R\'eciproquement, s'il existe un \el\ $\w$ hors de $<h'>$ qui conjugue deux \el s non triviaux de $<h'>$, alors soit
$\w \in <aha^{-1}>$, soit il existe un \el\ hors de $<h>$ conjuguant deux \el s non triviaux de $<h>$. Ainsi il suffit
pour en d\'ecider  de combiner l'algorithme du théorème \ref{gwp_gro}, \`a une r\'esolution  pour $h'$.  $\ \square$

\begin{thm}[{\it Malnormalité d'un famille finie de sous-groupe cycliques infinis}]
\label{lim_gro} Soit $G$ un groupe hyperboli\-que, et $h_1,h_2\in G$ des \el s d'ordre infini. On peut d\'ecider si  il
existe un \el\ non trivial  de $<h_1>$  conjugu\'e \`a un \el\ de $<h_2>$. En particulier, donnée une famille finie de
sous-groupes cycliques infinis de $G$, on peut décider si cette famille est malnormale dans $G$.

\end{thm}

\noindent \textsl{D\'emonstration.} On se donne $\l\geq 1,\e\geq
0$, et $C$ comme dans le lemme \ref{Creduit}.

Si $h_1$ et $h_2$ sont des mots $C$-r\'eduits, le lemme \ref{V3}
donne une constante calculable $L$, telle que si un \el\ non
trivial de $<h_1>$ est conjugu\'e \`a un \el\ de $h_2$, il existe
$n_1,n_2\in \Z_\ast$, tels que $|n_1|<L,|n_2|<L$, et $h^{n_1}$ est
conjugu\'e \`a $h^{n_2}$. On peut alors utiliser l'algorithme de
la conjugaison pour d\'ecider si c'est ou non le cas.

Si $h_1$ ou $h_2$ n'est pas $C$-r\'eduit. On applique le proc\'ed\'e de $C$-r\'eduction jumel\'ee, pour obtenir les
mots $C$-r\'eduits $H_1,H_2$, conjugu\'es respectifs de $h_1^N$ et $h_2^N$ pour un certain $N$. Clairement, un \el\ non
trivial de $<h_1>$ est conjugu\'e \`a un \el\ de $<h_2>$ si et seulement si un \el\ non trivial de $<H_1>$ est
conjugu\'e \`a un \el\ de $<H_2>$. On peut donc conclure. $\ \square$

\subsection{Centre, racines, normalisateur et centralisateur}

\begin{thm}[{\it Centralisateur et normalisateur d'un sous-groupe
cyclique infini}] \label{Z(h)} Soient $G$ un groupe hyperbolique et $h\in G$ un \el\ d'ordre infini. On peut
algorithmiquement déterminer le centralisateur $Z_G(h)$ et le normalisateur $N_G(h)$ de $<h>$. Plus précisemment, ils
contiennent tous deux $<h>$ comme sous-groupe d'indice fini et l'on peut algorithmiquement déterminer des représentants
de leurs classes à gauche modulo $<h>$.
\end{thm}

\noindent {\sl Démonstration.} Dans la suite on note indifféremment $\Lambda(h)$ pour le centralisateur et le
normalisateur de $<h>$, l'argument étant identique dans les deux cas en prenant soin de noter $e=1$ lorsque
$\Lambda(h)=Z_G(h)$ et $e=\pm 1$ lorsque $\Lambda(h)=N_G(h)$. On se donne une famille génératrice de $G$ pour laquelle
$G$ est $\delta$-hyperbolique, $\d$ et les constantes $C$, $\l$, et $\e$ provenant du lemme \ref{Creduit}.

a) Si $h$ est $C$-réduit. Clairement $<h>\subset \Lambda(h)\subset Stab(h^-,h^+)$ ; avec le lemme \ref{lemstab} il
existe une constante calculable $k=k(\d,\l,\e)$ tel que toute classe à gauche de $Stab(h^-,h^+) \mod <h>$ a un
représentant dans $Stab(h^-,h^+)$ de longueur au plus $k$. On utilise l'algorithme du mot pour décider si
$uhu^{-1}h^e=1$ pour chaque $u\in G$ dans la boule de rayon $k$. Celà permet de déterminer $u_1,\ldots , u_p \in G$ qui
avec $h$ forment une famille génératrice de $\Lambda(h)$. En utilisant une solution au problème du mot généralisé de
$<h>$ (théorème \ref{gwp_gro}) on détermine une famille $1,u_1,\ldots , u_q$ de représentants des classes à gauche de
$\Lambda(h) \mod <h>$.

b) Si $h$ n'est pas $C$-réduit on applique la procédure de $C$-réduction pour obtenir $h_1$ $C$-réduit, $n\in \N^*$, et
$a\in G$ tels que $h_1=ah^na^{-1}$. On détermine $\Lambda(h_1)$ avec a) puis $\Lambda(h^n)=a^{-1}\Lambda(h_1)a$. Celà
fournit des représentants $1,v_1,\ldots, v_r$ des classes à gauche de $\Lambda(h^n)\mod <h^n>$, et on en déduit avec le
théorème \ref{gwp_gro} des représentants $1,v_1,\ldots,v_s$ des classes à gauche de $\Lambda(h^n)\mod <h>$. Puisque
$<h>\subset\Lambda (h)\subset\Lambda(h^n)$ il suffit de décider avec l'algorithme du mot pour $i=1,\ldots ,s$ si
$v_ihv_ih^e=1$ pour obtenir les représentants $1,v_1,\ldots ,v_t$ de $\Lambda (h)\mod <h>$.$\ \square$

\begin{thm}[{\it Racines d'un \el\ d'ordre infini}]\label{root}
Soit $G$ un groupe hyperbolique et $h\in G$ un \el\ d'ordre infini. On peut algorithmiquement déterminer toutes les
racines de $h$ (leur nombre est fini, majoré par l'indice de $<h>$ dans $Z(h)$).
\end{thm}

\noindent {\sl Démonstration.} Toutes les racines de $h$ sont contenues dans $Z(h)$. On détermine des représentants $1,
a_1, \ldots , a_r$ des classes à gauche de $Z(h)\mod <h>$ avec l'algorithme du théorème précédent.  Pour tout
$i=1,\ldots, r$ il existe $p_i\in\N^*$ tel que $a_i^{n}\in <h>$ si et seulement si $p_i$ divise $n$ ; posons
$a_i^{p_i}=h^{q_i}$. Soit $r_i=h^{x_i}a_i\in Z(h)$ ; $r_i^n\in <h>$ si et seulement si $p_i$ divise $n$.
$$r_i^{p_i}=(h^{x_i}a_i)^{p_i}=h^{p_ix_i+q_i}$$
Ainsi $r_i$ est une racine de $h$ si et seulement si $p_ix_i+q_i=1$ c'est à dire si et seulement si $x_i=(1-q_i)/p_i$.
Ainsi pour chaque $i=1,\ldots ,r$, $h$ admet une racine $r_i$ si et seulement si $p_i$ divise $1-q_i$. On utilise les
algorithmes des théorèmes \ref{torsion} et \ref{gwp_gro} avec l'algorithme du mot pour déterminer chacun des $p_i$ et
$q_i$. Celà permet de déterminer toutes les racines de $h$.$\ \square$

\begin{thm}[{\it Maximalité d'un sous-groupe cyclique infini}]
Soient $G$ un groupe hyperbolique et $H$ un sous-groupe cyclique infini de $G$. On peut décider si $H$ est maximal dans
$G$.
\end{thm}

\noindent {\sl Démonstration.} On se donne un \el\ arbitraire de $H$ que l'on choisit non trivial à l'aide de
l'algorithme du mot. Il suffit d'appliquer l'algorithme du théorème \ref{root} pour déterminer toutes ses racines puis
pour chacune d'entre-elles l'algorithme de \ref{gwp_gro} pour décider si elle est dans $H$.$\ \square$

\begin{thm}[{\it Centralisateur d'un sous-groupe de type fini qui n'est pas de torsion}]
Soit $G$ un groupe hyperbolique et $K$ un sous-groupe de type fini qui n'est pas de torsion. On peut algorithmiquement
déterminer le centralisateur $Z_G(K)$ de $K$ dans $G$.
\end{thm}

\noindent {\sl Démonstration.} Soient $k_1,\ldots , k_n$ des générateurs de $K$ et $k\in K$ un \el\ d'ordre infini.
Clairement $Z_G(K)\subset Z_G(k)$. On utilise l'algorithme du théorème \ref{Z(h)} pour déterminer des représentants
$1,a_1,\ldots ,a_r$ de $Z_G(k)\mod <k>$. On utilise l'algo\-rithme du théorème \ref{2coset_gro} pour déterminer pour
chaque $i=1,\ldots ,r$ tous les entiers $n$ tels que $a_ik^na_i^{-1}=k^n$ ; celà permet de déterminer $Z_G(K)\cap <k>$.

Si $Z_G(K)\cap <k>=\{1\}$ alors $Z_G(K)$ est fini. Chaque classe de $Z_G(k)\mod <k>$ contient au plus un \el\ de
torsion que l'on détermine aisément en trouvant pour chaque $a_i$ les entiers $p_i$ et $q_i$ définis comme dans la
preuve du théorème \ref{root}. Il ne reste alors qu'à utiliser l'algorithme du mot pour décider pour chacun d'entre eux
s'il commute avec les générateurs de $K$.

Si $Z_G(K)\cap <K>=<k^n>$ ; on déduit de $a_1,\ldots ,a_r$ des représentants $1,b_1,\ldots, b_s$ des classes à gauche
de $Z_G(k)\mod <k^n>$. Etre dans $Z_G(K)$ ne dépend pas du représentant $\mod <k^n>$ considéré. Pour chaque $b_i$ on
décide avec l'algorithme du mot s'il commute avec les générateurs de $K$ ce qui permet de déterminer les représentants
$1,b_1, \ldots , b_t$ des classes à gauche de $Z_G(K)\mod <k^n>$.$\ \square$

\begin{cor}[{\it Déterminer le centre}]
Soit $G$ un groupe hyperbolique qui n'est pas de torsion. On peut algorithmiquement déterminer le centre de $G$.
\end{cor}

\vskip 0.4cm


\end{document}